\documentclass[12pt]{article}
\usepackage{amssymb,latexsym,pstcol,pst-plot,amsthm,amsmath,hyperref}

\newcommand{\ben}{\begin{enumerate}}
\newcommand{\een}{\end{enumerate}}
\newcommand{\ble}{\begin{lem}}
\newcommand{\ele}{\end{lem}}
\newcommand{\bth}{\begin{thm}}
\renewcommand{\eth}{\end{thm}}
\newcommand{\bpr}{\begin{prop}}
\newcommand{\epr}{\end{prop}}
\newcommand{\bco}{\begin{cor}}
\newcommand{\eco}{\end{cor}}
\newcommand{\bcon}{\begin{conj}}
\newcommand{\econ}{\end{conj}}
\newcommand{\bde}{\begin{defn}}
\newcommand{\ede}{\end{defn}}
\newcommand{\bex}{\begin{exa}}
\newcommand{\eex}{\end{exa}}
\newcommand{\barr}{\begin{array}}
\newcommand{\earr}{\end{array}}
\newcommand{\btab}{\begin{tabular}}
\newcommand{\etab}{\end{tabular}}
\newcommand{\beq}{\begin{equation}}
\newcommand{\eeq}{\end{equation}}
\newcommand{\bea}{\begin{eqnarray*}}
\newcommand{\eea}{\end{eqnarray*}}
\newcommand{\bce}{\begin{center}}
\newcommand{\ece}{\end{center}}
\newcommand{\bpi}{\begin{picture}}
\newcommand{\epi}{\end{picture}}
\newcommand{\bpp}{\begin{picture}}
\newcommand{\epp}{\end{picture}}
\newcommand{\bfi}{\begin{figure} \begin{center}}
\newcommand{\efi}{\end{center} \end{figure}}
\newcommand{\bprf}{\begin{proof}}
\newcommand{\eprf}{\end{proof}\medskip}
\newcommand{\capt}{\caption}
\newcommand{\bsl}{\begin{slide}{}}
\newcommand{\esl}{\end{slide}}
\newcommand{\bfr}{\begin{frame}}
\newcommand{\efr}{\end{frame}}

\newcommand{\prf}{{\bf Proof.\hspace{7pt}}}
\newcommand{\hqed}{\hfill \qed}
\newcommand{\hqedm}{\hfill \qed \medskip}

\newcommand{\hs}[1]{\hspace{#1}}
\newcommand{\hso}[1]{\hspace{-1pt}}

\newcommand{\emp}{\emptyset}

\newcommand{\sbe}{\subseteq}

\def\<{\langle}
\def\>{\rangle}

\newcommand{\ree}[1]{(\ref{#1})}

\newcommand{\ep}{\epsilon}

\newcommand{\la}{\lambda}

\newcommand{\cC}{{\cal C}}

\newcommand{\cL}{{\cal L}}

\newcommand{\Fib}[1]{\left\{ #1 \right\}}
\newcommand{\bin}[2]{\left\{ #1 \atop #2 \right\}}
\newcommand{\Luc}[1]{\langle #1 \rangle}

\newcommand{\Gbb}{\put(10,10){\circle*{3}}}

\newcommand{\Gdb}{\put(30,10){\circle*{3}}}

\newcommand{\Gfb}{\put(50,10){\circle*{3}}}

\newcommand{\Gaaea}{\put(0,0){\line(1,0){40}}}

\newcommand{\Gaaga}{\put(0,0){\line(1,0){60}}}

\newcommand{\Gaaac}{\put(0,0){\line(0,1){20}}}

\newcommand{\Gaaaf}{\put(0,0){\line(0,1){50}}}

\newcommand{\Gbabf}{\put(10,0){\line(0,1){50}}}

\newcommand{\Gcacc}{\put(20,0){\line(0,1){20}}}

\newcommand{\Gcacf}{\put(20,0){\line(0,1){50}}}

\newcommand{\Gdadf}{\put(30,0){\line(0,1){50}}}

\newcommand{\Geaec}{\put(40,0){\line(0,1){20}}}

\newcommand{\Geaef}{\put(40,0){\line(0,1){50}}}

\newcommand{\Ggagc}{\put(60,0){\line(0,1){20}}}

\newcommand{\Gabbb}{\put(0,10){\line(1,0){10}}}

\newcommand{\Gabeb}{\put(0,10){\line(1,0){40}}}

\newcommand{\Gbbdb}{\put(10,10){\line(1,0){20}}}

\newcommand{\Gdbfb}{\put(30,10){\line(1,0){20}}}

\newcommand{\Gfbgb}{\put(50,10){\line(1,0){10}}}

\newcommand{\Gacec}{\put(0,20){\line(1,0){40}}}

\newcommand{\Gacgc}{\put(0,20){\line(1,0){60}}}

\newcommand{\Gaded}{\put(0,30){\line(1,0){40}}}

\newcommand{\Gaeee}{\put(0,40){\line(1,0){40}}}

\newcommand{\Gafef}{\put(0,50){\line(1,0){40}}}

\newtheorem{thm}{Theorem}[section]
\newtheorem{prop}[thm]{Proposition}
\newtheorem{cor}[thm]{Corollary}
\newtheorem{lem}[thm]{Lemma}
\newtheorem{conj}[thm]{Conjecture}
\newtheorem{exa}[thm]{Example}

\begin{document}
\pagestyle{plain}

\title{Combinatorial interpretations of binomial coefficient analogues related
  to Lucas sequences
}
\author{
Bruce E. Sagan\footnote{Work partially done while a Program Officer at NSF.  The views
    expressed are not necessarily those of the NSF.}\\[-5pt]
\small Department of Mathematics, Michigan State University,\\[-5pt]
\small East Lansing, MI 48824-1027, USA, \texttt{sagan@math.msu.edu}\\
and\\
Carla D. Savage\footnote{Partially supported by NSA grant H98230-08-1-0072}\\[-5pt]
\small Department of Computer Science, North Carolina State University\\[-5pt]
\small Raleigh, NC 27695-8206, USA, \texttt{savage@cayley.csc.ncsu.edu}
}

\date{\today\\[10pt]
	\begin{flushleft}
	\small Key Words: combinatorial interpretation, binomial coefficient,
        fibonomial coefficient, integer partition, Lucas sequence
	                                       \\[5pt]
	\small AMS subject classification (2000):
	Primary 05A10;
	Secondary 05A17.
	\end{flushleft}}

\maketitle

\begin{abstract}
Let $s$ and $t$ be variables.  Define polynomials $\Fib{n}$ in $s,t$ by
$\Fib{0}=0$, $\Fib{1}=1$,
and $\Fib{n}=s\Fib{n-1}+t\Fib{n-2}$ for $n\ge2$.  If $s,t$ are integers then the corresponding
sequence of integers is called a {\it Lucas sequence\/}.  Define 
an analogue of the binomial coefficients by
$$
\bin{n}{k}=\frac{\Fib{n}!}{\Fib{k}!\Fib{n-k}!}
$$
where $\Fib{n}!=\Fib{1}\Fib{2}\cdots\Fib{n}$.  It is easy to see that $\bin{n}{k}$ is a
polynomial in $s$ and $t$.  The purpose of this note is to give two
combinatorial interpretations for this polynomial in terms of statistics on
integer partitions inside a $k\times (n-k)$ rectangle.  When $s=t=1$ we obtain combinatorial
interpretations of the fibonomial coefficients which are simpler than any that
have previously appeared in the literature.
\end{abstract}

\section{Introduction}

Given variables $s,t$ we define the corresponding sequence of {\it Lucas polynomials\/}, $\Fib{n}$, by  $\Fib{0}=0$, $\Fib{1}=1$, and for $n\ge2$:
\beq
\label{Fib}
\Fib{n}=s\Fib{n-1}+t\Fib{n-2}.
\eeq
When $s,t$ are integers, the corresponding integer sequence is called a
{\it Lucas sequence\/} \cite{luc:tfn1,luc:tfn2,luc:tfn3}.  These sequences  have many interesting
number-theoretic and combinatorial properties.

Define the {\it lucasnomials\/}  for $0\le k\le n$ by
\beq
\label{bineq}
\bin{n}{k}=\frac{\Fib{n}!}{\Fib{k}!\Fib{n-k}!}
\eeq
where $\Fib{n}!=\Fib{1}\Fib{2}\cdots\Fib{n}$.  It is not hard to show that $\bin{n}{k}$ is a
polynomial in $s$ and $t$. (This follows from Proposition~\ref{binpr} below).  The purpose of
this note is to give two simple combinatorial interpretations of
the lucasnomials.  They  are based on statistics associated with
integer partitions $\la$ inside a $k\times (n-k)$ rectangle.  More specifically, we will
show that $\bin{n}{k}$ is the generating function for
certain tilings of such $\la$  and their complements  with dominos and
monominos.

Various specializations of the parameters $s$ and $t$ are of interest.
When $s=t=1$, $\Fib{n}$ becomes the $n$th Fibonacci number, and
$\bin{n}{k}$ is known as a fibonomial coefficient.  Gessel and Viennot~\cite{gv:bdp}
gave an interpretation of the fibonomials in terms of nonintersecting lattice paths
and asked for a simpler one.  Benjamin and Plott~\cite{bp:caf}
gave another interpretation in terms of tilings, but it is not as
straight-forward as ours.   When $s=\ell$ and $t=-1$ we get two new
interpretations of the $\ell$-nomial
coefficients of Loehr and Savage~\cite{ls:gcb}.  This case is of interest because
of its connection with the Lecture Hall Partition Theorem introduced by
Bousquet-M\'elou and (Kimmo) Eriksson~\cite{bme:lhp1,bme:lhp2,bme:rlh}.  Finally, letting $s=q+1$
and $t=-q$ one gets new interpretations for the classical $q$-binomial
coefficients.  For more information about these important polynomials, see the
text of Andrews~\cite{and:tp}.

\section{Recursions}

In this section we will present the recurrence relations we will need for
our combinatorial interpretations.  To obtain these results, we will use a
tiling model that is often useful when dealing with Lucas sequences.  See the
book of Benjamin and Quinn~\cite{bq:prc} for more details.  Suppose we have $n$ squares
arranged in a $1\times n$ rectangle.  We number the squares $1,\ldots,n$ from
left to right and also number the vertical edges of the squares $0,\ldots,n$
left to right.  A {\it linear tiling\/}, $T$, is a covering of the rectangle
with disjoint dominos (covering two squares) and monominos (covering one
square).  Let $\cL_n$ be the set of all such $T$.  The tilings in $\cL_3$ are
drawn in Figure~\ref{L3} where a dot in a square represents a monomino while
two dots connected by a horizontal line represent a domino.

\thicklines
\setlength{\unitlength}{1pt}
\bfi
$$
\barr{ccccccc}
\cL_3&:&
\bpi(60,20)(0,0)
\Gaaga \Gacgc
\Gaaac \Gcacc \Geaec \Ggagc
\Gbb \Gdb \Gfb
\epi
&&
\bpi(60,20)(0,0)
\Gaaga \Gacgc
\Gaaac \Gcacc \Geaec \Ggagc
\Gbb \Gdb \Gfb \Gbbdb
\epi
&&
\bpi(60,20)(0,0)
\Gaaga \Gacgc
\Gaaac \Gcacc \Geaec \Ggagc
\Gbb \Gdb \Gfb \Gdbfb
\epi
\\[5pt]
\Fib{4}&=&s^3&+&st&+&st
\earr
$$
\capt{The tilings in $\cL_3$ and corresponding weights}
\label{L3}
\efi

Let the weight of a tiling be $w(T)=s^m t^d$ where $m$ and $d$ are the number of
monominos  and dominos in $T$, respectively.  We will use the same weight
for all other types of tilings considered below.  Since the last tile in any
tiling must be a monomino or domino, the initial conditions and
recursion~\ree{Fib} immediately give
$$
\Fib{n+1}=\sum_{T\in \cL_n} w(T).
$$
See Figure~\ref{L3} for an illustration of the case $n=3$.

For our second combinatorial interpretation, we will need another sequence of
polynomials closely related to the $\Fib{n}$.  Define 
$\Luc{n}$ using recursion~\ree{Fib} but with the initial
conditions $\Luc{0}=2$ and $\Luc{1}=s$.  If $s=t=1$ then $\Luc{n}$ is the
$n$th Lucas number.  A combinatorial interpretation for these polynomials is
obtained via another type of tiling.  In a {\it circular tiling\/} of the
$1\times n$ rectangle, the edges labeled $0$ and and $n$ are identified so
that it is possible to have a domino crossing this edge and covering the
first and last squares.  Such a domino, if it exists, will be called the {\it circular domino\/} of the tiling.   Let $\cC_n$ be the set of circular tilings of a $1\times n$ rectangle.  So $\cL_n\sbe\cC_n$ is the subset of all circular tilings with no circular domino.  For example, $\cC_3$ consists of the tilings in $\cL_3$ displayed previously together with
$$
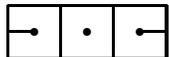
\bpi(60,20)(0,0)
\Gaaga \Gacgc
\Gaaac \Gcacc \Geaec \Ggagc
\Gbb \Gdb \Gfb \Gabbb \Gfbgb
\epi
.
$$

Now for $n\ge1$ we have
\beq
\label{Luc}
\Luc{n}=\sum_{T\in\cC_n} w(T).
\eeq
Indeed, to show that the sum satisfies~\ree{Fib} first note that we already
have a weight-preserving bijection for the linear tilings involved.  And if
$T\in\cC_n$ has a circular edge, then removal of the tile covering square
$n-1$ will take care of the remainder.  In order to make~\ree{Luc} also hold
for $\cC_0$, we give the empty tiling $\ep$ of the $1\times 0$ box weight
$w(\ep)=2$.  Bear in mind that $\ep$ considered as an element of $\cL_0$ still
has $w(\ep)=1$.  Context will always make it clear which weight we are using. 

We start with two recursions for the Lucas polynomials.  These are well known
for Lucas sequences; see~\cite[p.\ 38, Identity 73]{bq:prc} for~\ree{fib1}
and~\cite[p.\ 46, Identity 94]{bq:prc} or~\cite[p.\ 201, Equation 49]{luc:tfn2} for~\ree{fib2}.  
Also, the proofs in the integer case generalize
to variable $s$ and $t$ without difficulty.  But we will provide a
demonstration for completeness and to emphasize the simplicity of the
combinatorics involved. 
\ble
\label{m+n}
For $m\ge1$ and $n\ge0$ we have
\beq
\label{fib1}
\Fib{m+n} = \Fib{n+1}\Fib{m}+t\Fib{m-1}\Fib{n}.
\eeq
For $m,n\ge0$ we have
\beq
\label{fib2}
\Fib{m+n}=\frac{\Luc{n}}{2} \Fib{m} + \frac{\Luc{m}}{2} \Fib{n}.
\eeq
\ele
\prf
For the first identity, the left-hand side is the generating function for $T\in\cL_{m+n-1}$.  The second and first terms on the right correspond to those tilings which do or do not have a domino crossing the edge labeled $n$, respectively.

Multiply the second equation by 2 and consider two copies of $\cL_{m+n-1}$.  In each tiling in the first copy distinguish the edge labeled $m-1$, and do the same for the edge labeled $m$ in the second copy.  The set of tilings in both copies where a domino does not cross the distinguished edge accounts for the terms corresponding to linear pairs in $\Luc{n}\Fib{m}+\Luc{m}\Fib{n}$.  If a domino crosses the distinguished edge $m-1$ in a tiling $T$, then consider the restriction $T'$ of $T$ to the first $m$ squares as a circular tiling by shifting it so that the domino between squares $m-1$ and $m$ becomes the circular edge.  Also let $T''$ be the restriction of $T$ to the last $n-1$ squares considered as a linear tiling.  The pairs $(T',T'')$ account for the remaining terms in $\Luc{m}\Fib{n}$.  A similar bijection using the tilings with a domino crossing the distinguished edge $m$ accounts for the rest of the terms in $\Luc{n}\Fib{m}$, completing the proof.\hqedm

We can use the recurrence relations in the previous lemma to produce recursions for the lucasnomials.

\bpr
\label{binpr}
For $m,n\ge1$ we have
\bea
\bin{m+n}{m} &=& \Fib{n+1}\bin{m+n-1}{m-1} + t\Fib{m-1}\bin{m+n-1}{n-1}\\[10pt]
&=&\frac{\Luc{n}}{2} \bin{m+n-1}{m-1} + \frac{\Luc{m}}{2} \bin{m+n-1}{n-1}.
\eea
\epr
\prf
Given $m,n$ and  any polynomials $p$ and $q$ such that $\Fib{m+n}= p \Fib{m}+
q \Fib{n}$, straightforward algebraic manipulation of the definition of
lucasnomials yields 
$$
\bin{m+n}{m} = p\bin{m+n-1}{m-1} + q\bin{m+n-1}{n-1}.
$$
Combining this observation with Lemma~\ref{m+n}, we are done. 
\hqedm

\section{The combinatorial interpretations}

Our combinatorial interpretations of $\bin{m+n}{m}$ will involve integer partitions.  A {\it partition\/} is a weakly decreasing sequence $\la=(\la_1,\la_2,\ldots,\la_m)$ of nonnegative integers.  The $\la_i$, $1\le i\le m$, are called {\it parts\/} and note that we are allowing zero as a part.  The {\it Ferrers diagram\/} of $\la$, also denoted $\la$, is an array of $m$ left-justified rows of boxes with $\la_i$ boxes in row $i$.
We say that $\la$ is {\it contained in an $m\times n$ rectangle\/}, written $\la\sbe m\times n$, if it has $m$ parts and each part is at most $n$.  In this case, $\la$ determines another partition $\la^*$ whose parts are the lengths of the
columns of the complement of $\la$ in $m\times n$.  The first diagram in Figure~\ref{la} shows $\la=(3,2,2,0,0)$ contained in a $5\times 4$ rectangle with complement $\la^*=(5,4,2,2)$.

\thicklines
\setlength{\unitlength}{2pt}
\bfi
$$
\barr{cc}
\bpi(60,50)(0,0)
\Gaaea \Gabeb \Gacec \Gaded \Gaeee \Gafef
\Gaaaf \Gbabf \Gcacf \Gdadf \Geaef
\multiput(0,0)(.1,0){11}{\line(0,1){20}}
\multiput(0,20)(0,-.1){11}{\line(1,0){21}}
\multiput(20,20)(.1,0){11}{\line(0,1){20}}
\multiput(20,40)(0,-.1){11}{\line(1,0){11}}
\multiput(30,40)(.1,0){11}{\line(0,1){10}}
\multiput(30,50)(0,-.1){11}{\line(1,0){10}}
\epi
&
\bpi(60,50)(0,0)
\Gaaea \Gabeb \Gacec \Gaded \Gaeee \Gafef
\Gaaaf \Gbabf \Gcacf \Gdadf \Geaef
\multiput(0,0)(.1,0){11}{\line(0,1){20}}
\multiput(0,20)(0,-.1){11}{\line(1,0){21}}
\multiput(20,20)(.1,0){11}{\line(0,1){20}}
\multiput(20,40)(0,-.1){11}{\line(1,0){11}}
\multiput(30,40)(.1,0){11}{\line(0,1){10}}
\multiput(30,50)(0,-.1){11}{\line(1,0){10}}
\multiput(5,5)(0,10){5}{\circle*{1.5}}
\multiput(15,5)(0,10){5}{\circle*{1.5}}
\multiput(25,5)(0,10){5}{\circle*{1.5}}
\multiput(35,5)(0,10){5}{\circle*{1.5}}
\multiput(5,5)(10,0){4}{\line(0,1){10}}
\put(5,35){\line(1,0){10}}
\put(15,45){\line(1,0){10}}
\put(35,35){\line(0,1){10}}
\epi
\earr
$$
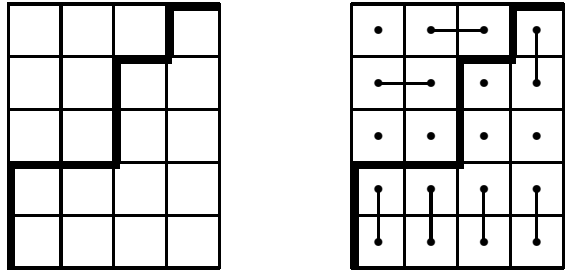
\capt{A partition $\la$ contained in $5\times 4$ and a tiling}
\label{la}
\efi

A {\it linear tiling of $\la$\/} is a covering of its Ferrers diagram with
disjoint dominos and monominos obtained by linearly tiling each $\la_i$.  The
set of such tilings is denoted $\cL_\la$.  Note that if $\la\sbe m\times n$
then  $T\in\cL_\la$ gives a tiling of each of its rows, while $T\in
\cL_{\la^*}$ gives a tiling of each of its columns.  We will also need
$\cL'_n$ which is the set of all tilings in $\cL_n$ which do not begin with a
monomino.  This is equivalent to beginning with a domino if $n\ge2$, and for
$n<2$ yields $\cL'_0=\{\ep\}$ and $\cL'_1=\emp$.  We define $\cL'_\la$ similarly.  The second diagram in
Figure~\ref{la} shows a tiling in $\cL_\la\times\cL'_{\la^*}$.  In a {\it
  circular tiling of $\la$\/} we use circular tilings on each $\la_i$.  So  if
$\la_i=0$ then it will get the empty tiling which has weight 2.  The notation
$\cC_\la$ is self-explanatory.  If one views the tiling $T$ in Figure~\ref{la}
as an element of $\cL_\la\times\cL'_{\la^*}$ then it has weight $w(T)=s^6
t^7$.  But as an element of $\cC_\la\times\cC_{\la^*}$ it has $w(T)=4 s^6
t^7$.  As usual, context will clarify which weight to use.  We are now in a
position to state and prove our two combinatorial interpretations for the
lucasnomials.  The first has the nice property that it is multiplicity free.
The second is pleasing because it displays the natural symmetry of
$\bin{m+n}{m}$. 

\bth
For $m,n\ge0$ we have
\beq
\label{bineq1}
\bin{m+n}{m}=\sum_{\la\sbe m\times n}\hs{5pt} \sum_{T\in\cL_\la\times\cL'_{\la^*}} w(T),
\eeq
and
\beq
\label{bineq2}
2^{m+n}\bin{m+n}{m}=
\sum_{\la\sbe m\times n} \hs{5pt} \sum_{T\in\cC_\la\times\cC_{\la^*}} w(T).
\eeq
\eth
\prf
We will show that the right-hand side of~\ree{bineq1} satisfies the first recursion in Proposition~\ref{binpr} as the initial conditions are easy to verify.  Given
$\la\sbe m\times n$ there are two cases.  If $\la_1=n$ then the generating function for tilings of the first row of $\la$ is $\Fib{n+1}$, and $\bin{m+n-1}{m-1}$ counts the ways to fill the rest of the rectangle.  If $\la_1<n$ then $\la^*_1=m$.  The generating function for $\cL'_m$ is $t\Fib{m-1}$ and $\bin{m+n-1}{n-1}$ takes care of the rest.

Proving that both sides of~\ree{bineq2} is similar using the fact that if we let $f(m,n)=2^{m+n}\bin{m+n}{m}$ then, by the second recursion in Proposition~\ref{binpr}, we have $f(m,n)=\Luc{n}f(m-1,n)+\Luc{m}f(m,n-1)$.  This completes the proof.\hqed

\end{document}